\documentclass[11pt]{amsart}
\usepackage{amsfonts, amsthm, amssymb, amsmath}
\usepackage[margin=1.1in]{geometry}
\usepackage{enumitem}

\newtheorem{theorem}{Theorem}[section]
\newtheorem{proposition}[theorem]{Proposition}
\newtheorem{lemma}[theorem]{Lemma}
\newtheorem{corollary}[theorem]{Corollary}
\theoremstyle{definition}

\newtheorem{question}[theorem]{Question}
\theoremstyle{remark}
\newtheorem{remark}[theorem]{Remark}

\newcommand{\NFS}{\mathrm{NFS}}
\newcommand{\MA}{\mathrm{MA}(\omega_1)}
\newcommand{\PFA}{\mathrm{PFA}}
\newcommand{\cl}[1]{\overline{#1}}
\newcommand{\R}{\mathbb R}
\newcommand{\N}{\mathbb N}
\newcommand{\es}{\varnothing}
\DeclareMathOperator{\hL}{hL}
\DeclareMathOperator{\hd}{hd}

\begin{document}

\title[HL and HS subspaces of function spaces]{Hereditarily Lindel\"of and hereditarily separable\\ subspaces of function spaces}

\author{Evgenii Reznichenko}
\address{Department of General Topology and Geometry, Faculty of Mechanics and Mathematics,\newline
Lomonosov Moscow State University, Leninskie Gory 1, Moscow, 119991 Russia}
\address{Shenzhen MSU-BIT University, Shenzhen, China}
\email{erezn@inbox.ru}


\begin{abstract}
A space is ESLC if it contains neither an $S$-space nor an $L$-space, and $\NFS$ is the class of spaces containing no
free sequences of length $\omega_1$. We prove that under $\MA$, if all finite powers of $X$ belong to $\NFS$, then
$C_p(X)^\omega$ is ESLC; thus, for subspaces of $C_p(X)$, hereditary Lindel\"ofness, hereditary separability and
countable spread coincide. This is derived from a dual form of a theorem of Todor\v{c}evi\'c. A Lindel\"of
$\Sigma$-space belongs to $\NFS$ if and only if its compact subspaces have countable tightness. Hence $C_p(X)$ contains
no $L$-spaces if $X$ is a Lindel\"of $\Sigma$-space of countable tightness and $\MA$ holds, or if $X$ is any
Lindel\"of $\Sigma$-space and $\PFA$ holds. Under $\MA$, all iterated function spaces of Corson compact spaces are
ESLC, and in $\mathrm{ZFC}$ the space $C_p(X,\{0,1\})$ is ESLC for every zero-dimensional Corson compact space $X$.
\end{abstract}

\maketitle

\section{Introduction}

All spaces are assumed to be Tychonoff, and $C_p(X)$ denotes the space of continuous real-valued functions on $X$
with the topology of pointwise convergence. Following \cite{Arh}, a space is \emph{Lindel\"of $\Sigma$} if it is a
continuous image of a Lindel\"of $p$-space, i.e., of a perfect preimage of a separable metrizable space; equivalently,
Lindel\"of $\Sigma$-spaces form the smallest class containing all compact and all separable metrizable spaces and
closed under countable products, closed subspaces and continuous images.

An \emph{$S$-space} is a hereditarily separable space which is not hereditarily Lindel\"of, and an \emph{$L$-space} is
a hereditarily Lindel\"of space which is not hereditarily separable. Under $\PFA$ there are no $S$-spaces \cite{T89},
while $L$-spaces exist in $\mathrm{ZFC}$ \cite{Moore}. Following Hart and Kunen \cite{HK}, we call a space \emph{ESLC}
if it contains neither an $S$-space nor an $L$-space. A space $X$ is ESLC if and only if for every $Y\subset X$ the
following conditions are equivalent: $Y$ is hereditarily separable; $Y$ is hereditarily Lindel\"of; $Y$ contains no
uncountable discrete subspaces \cite[\S1]{HK}.\footnote{If $Y$ is not hereditarily Lindel\"of, it contains a subset
$Y'$ which is right-separated in type $\omega_1$; if $Y'$ is not hereditarily separable either, it contains a subset
$D$ which is left-separated in type $\omega_1$ with respect to another well-ordering, and by the partition relation
$\omega_1\to(\omega_1,\omega)^2$ of Dushnik and Miller we may assume that the two well-orderings agree on $D$; then
$D$ is discrete.}

A sequence $(x_\alpha)_{\alpha<\omega_1}$ in a space $X$ is \emph{free} if
\[
\cl{\{x_\alpha:\alpha<\beta\}}\cap\cl{\{x_\alpha:\alpha\ge\beta\}}=\es
\qquad\text{for every } \beta<\omega_1 ,
\]
and $\NFS$ is the class of spaces containing no free sequences of length $\omega_1$. Every Lindel\"of space of
countable tightness belongs to $\NFS$, a compact space belongs to $\NFS$ if and only if it has countable tightness, and
$\NFS$ is preserved by closed subspaces and continuous images (Lemma~\ref{l:NFS}). By a theorem of Todor\v{c}evi\'c
\cite[Theorem 7.10]{T89} (see also \cite[Theorem 13]{T93}), under $\MA$ every space all of whose finite powers belong to
$\NFS$ is ESLC. Since the finite powers of such a space satisfy the same assumption, they are ESLC as well, and so is
its countable power, because ESLC is quasi-productive: a countable product is ESLC whenever all its finite subproducts
are \cite[Lemma 3.2]{HK}. For function spaces, results of this kind go back to the theorem of the author that under
$\MA$ every separable compact space $K$ with Lindel\"of $C_p(K)$ is metrizable (see \cite{Rez92} and
\cite[Ch.~IV, \S8]{Arh}), and to the related theorems of Okunev~\cite{Ok} and of Okunev and the author~\cite{OR}.

Our main result is the following theorem.

\begin{theorem}\label{t:main}
$(\MA)$ Let $X$ be a space such that all finite powers of $X$ belong to $\NFS$. Then $C_p(X)^\omega$ is ESLC. In
particular, for every $Y\subset C_p(X)$ the following conditions are equivalent:
\begin{enumerate}[label=\textup{(\arabic*)}]
\item $Y$ contains no uncountable discrete subspaces;
\item $Y$ is hereditarily Lindel\"of;
\item $Y$ is hereditarily separable.
\end{enumerate}
\end{theorem}

The assumption $\MA$ cannot be omitted: under $\mathfrak b=\omega_1$, in particular under $\mathrm{CH}$,
Theorems~\ref{t:main}--\ref{t:B} and Corollaries~\ref{c:Sigma} and~\ref{c:L}(1) fail (Section~\ref{s:noMA}).
Theorem~\ref{t:main} is a dual form of the following theorem of Todor\v{c}evi\'c, in which the roles of $X$ and
$C_p(X)$ are interchanged.

\begin{theorem}[{\cite[Theorem 19]{T93}}]\label{t:T19}
$(\MA)$ If $C_p(X)\in\NFS$, then $X$ is ESLC.
\end{theorem}

A family $\Phi\subset C_p(X)$ \emph{generates the topology} of $X$ if the sets $\varphi^{-1}(V)$, where
$\varphi\in\Phi$ and $V\subset\R$ is open, form a subbase of $X$. The following theorem is the main technical result
of the paper.

\begin{theorem}\label{t:A}
$(\MA)$ Let $\Phi\subset C_p(X)$ generate the topology of $X$, and suppose that $\Phi^n\in\NFS$ for all $n\in\N$.
Then $X$ is ESLC.
\end{theorem}

Theorem~\ref{t:main} is obtained by applying Theorem~\ref{t:A} to the space $C_p(X)$ and the family
$\{\hat x:x\in X\}\subset C_p(C_p(X))$ of evaluations, which is a continuous image of $X$. Conversely, if $\Phi$
generates the topology of $X$, then $X$ embeds in $C_p(\Phi)$, so Theorems~\ref{t:main} and~\ref{t:A} are
equivalent.

Theorem~\ref{t:A} with $\Phi=C_p(X)$ does not imply Theorem~\ref{t:T19}, since it requires all finite powers of
$C_p(X)$ to belong to $\NFS$. Todor\v{c}evi\'c's proof of Theorem~\ref{t:T19} uses sums of functions, which yields a
free sequence in $C_p(X)$ itself rather than in its finite powers. The same argument gives the following version of
Theorem~\ref{t:A}; applied to $\Phi=C_p^+(X)=\{f\in C_p(X):f\ge0\}$, which is closed in $C_p(X)$, it yields
Theorem~\ref{t:T19}. A common generalization, which also covers families closed under finite maxima, is
Theorem~\ref{t:Bgen}.

\begin{theorem}\label{t:B}
$(\MA)$ Let $\Phi\subset C_p(X)$ be a family of nonnegative functions closed under finite sums such that for every
closed $F\subset X$ and $x\in X\setminus F$ there is $f\in\Phi$ with $f(x)\ge1$ and $f(F)\subset\{0\}$. If
$\Phi\in\NFS$, then $X$ is ESLC.
\end{theorem}

For Lindel\"of $\Sigma$-spaces, the property $\NFS$ is determined by compact subspaces and is preserved by countable
powers; the following theorem is proved in $\mathrm{ZFC}$.

\begin{theorem}\label{t:Sigma}
For a Lindel\"of $\Sigma$-space $X$ the following conditions are equivalent:
\begin{enumerate}[label=\textup{(\arabic*)}]
\item $X\in\NFS$;
\item all compact subspaces of $X$ have countable tightness;
\item $X$ is a continuous image of a Lindel\"of $p$-space of countable tightness;
\item $X^\omega\in\NFS$.
\end{enumerate}
\end{theorem}

We recall the following known facts.

\begin{proposition}\label{p:known}
Let $X$ be a Lindel\"of $\Sigma$-space and $Y\subset C_p(X)$.
\begin{enumerate}[label=\textup{(\arabic*)}]
\item If $Y$ is separable, then $Y$ has a countable network.
\item If $Y$ contains no uncountable closed discrete subspaces, then $Y$ is Lindel\"of.
\item If $Y$ contains no uncountable discrete subspaces, then $Y$ is hereditarily Lindel\"of.
\end{enumerate}
\end{proposition}

Item (1) holds because $X$ is $\aleph_0$-stable, so that $C_p(X)$ is $\aleph_0$-monolithic \cite[Sect.~II.6]{Arh}
(see Section~\ref{s:prel}). Item (2) is Baturov's theorem \cite{Bat}, \cite[Theorem III.6.1]{Arh}, and (3) follows
from (2).

\begin{corollary}\label{c:Sigma}
$(\MA)$ Let $X$ be a Lindel\"of $\Sigma$-space in $\NFS$. Then every subspace of $C_p(X)$ containing no uncountable
discrete subspaces has a countable network.
\end{corollary}

\begin{corollary}\label{c:L}
Let $X$ be a Lindel\"of $\Sigma$-space, and suppose that either
\begin{enumerate}[label=\textup{(\arabic*)}]
\item $\MA$ holds and $t(X)\le\omega$, or
\item $\PFA$ holds.
\end{enumerate}
Then every subspace of $C_p(X)$ containing no uncountable discrete subspaces \textup(in particular, every hereditarily
Lindel\"of subspace\textup) has a countable network. In particular, $C_p(X)$ is ESLC.
\end{corollary}

For compact $X$, Corollary~\ref{c:L}(1) is essentially \cite[Theorem 1.5]{OR}.

In Section~\ref{s:Corson} we consider Corson compact spaces $X$. Under $\MA$, all finite and countable powers of the
iterated function spaces $C_{p,n}(X)$ are ESLC (Theorem~\ref{t:Corson}). In $\mathrm{ZFC}$, no $C_{p,n}(X)$ contains an
$S$-space, and if $X$ is zero-dimensional, then $C_p(X,\{0,1\})$ is ESLC; in fact, all its hereditarily Lindel\"of
subspaces are countable (Proposition~\ref{p:CpD}). On the other hand, consistently there is a Corson compact $L$-space
(Proposition~\ref{p:CorsonL}).

The paper is organized as follows. Section~\ref{s:prel} contains preliminaries. Theorems~\ref{t:A}, \ref{t:B} and
\ref{t:Bgen} are proved in Section~\ref{s:A}, and Theorem~\ref{t:main} in Section~\ref{s:main}. Lindel\"of
$\Sigma$-spaces are considered in Section~\ref{s:Sigma}, and the corollaries, including a theorem of Okunev, are proved
in Section~\ref{s:cor}. Section~\ref{s:Corson} deals with Corson and Gul'ko compact spaces, Section~\ref{s:noMA}
contains counterexamples without $\MA$, and Section~\ref{s:q} contains open questions. The results of this paper were
announced at the 2014 International Conference on Topology and its Applications (Nafpaktos, Greece).

\section{Preliminaries}\label{s:prel}

By $t(X)$, $s(X)$, $\hd(X)$ and $\hL(X)$ we denote the tightness, the spread (the supremum of cardinalities of
discrete subspaces), the hereditary density and the hereditary Lindel\"of number of $X$. A space $X$ is
\emph{$\aleph_0$-monolithic} if the closure of every countable subset of $X$ has a countable network, and
\emph{$\aleph_0$-stable} if every continuous image of $X$ which condenses onto a second countable space has a countable
network. By Arhangel'skii's duality theorems \cite[Theorems II.6.8 and II.6.9]{Arh}, $C_p(X)$ is $\aleph_0$-monolithic
if and only if $X$ is $\aleph_0$-stable, and $C_p(X)$ is $\aleph_0$-stable if and only if $X$ is $\aleph_0$-monolithic.
Compact spaces and, more generally, Lindel\"of $\Sigma$-spaces are $\aleph_0$-stable \cite[Sect.~II.6]{Arh}.

For a tuple $\bar f=(f_1,\dots,f_k)$ of functions and a point $x$ we write $\bar f(x)=(f_1(x),\dots,f_k(x))\in\R^k$.
For $A,B\subset\omega_1$ we write $A<B$ if $\alpha<\beta$ for all $\alpha\in A$ and $\beta\in B$.

\subsection{Free sequences}

We say that a sequence $(x_\alpha)_{\alpha<\omega_1}$ \emph{splits at} $\beta<\omega_1$ if
\[
\cl{\{x_\alpha:\alpha<\beta\}}\cap\cl{\{x_\alpha:\alpha\ge\beta\}}=\es.
\]
Thus, a sequence is free if it splits at every $\beta<\omega_1$.

\begin{lemma}\label{l:succ}
If a sequence $(x_\alpha)_{\alpha<\omega_1}$ splits at every successor ordinal, i.e.,
\[
\cl{\{x_\alpha:\alpha\le\gamma\}}\cap\cl{\{x_\alpha:\alpha>\gamma\}}=\es\qquad\text{for all }\gamma<\omega_1,
\]
then the sequence $y_\xi=x_{\xi+1}$, $\xi<\omega_1$, is free.
\end{lemma}

\begin{proof}
Let $\zeta<\omega_1$. Then $\{y_\xi:\xi<\zeta\}\subset\{x_\alpha:\alpha\le\zeta\}$ and
$\{y_\xi:\xi\ge\zeta\}\subset\{x_\alpha:\alpha>\zeta\}$, and it remains to apply the hypothesis with $\gamma=\zeta$.
\end{proof}

\begin{lemma}\label{l:NFS}
\begin{enumerate}[label=\textup{(\alph*)}]
\item Closed subspaces and continuous images of spaces in $\NFS$ are in $\NFS$.
\item Every Lindel\"of space of countable tightness is in $\NFS$.
\item A compact space $K$ is in $\NFS$ if and only if $t(K)\le\omega$.
\item Every space of countable spread is in $\NFS$.
\end{enumerate}
\end{lemma}

\begin{proof}
(a) Closures in a closed subspace coincide with closures in the whole space. If $g\colon X\to Y$ is a continuous
surjection and $(y_\alpha)$ is free in $Y$, choose $x_\alpha\in g^{-1}(y_\alpha)$; since
$\cl{\{x_\alpha:\alpha\in A\}}\subset g^{-1}\bigl(\cl{\{y_\alpha:\alpha\in A\}}\bigr)$, the sequence $(x_\alpha)$ is free in $X$.

(b) Let $(x_\alpha)$ be free in a Lindel\"of space $X$ with $t(X)\le\omega$. The points $x_\alpha$ are distinct, so
the set $\{x_\alpha:\alpha<\omega_1\}$ has a complete accumulation point $p$. Then $p\in\cl{\{x_\alpha:\alpha\ge\gamma\}}$
for all $\gamma$. By countable tightness, $p\in\cl C$ for a countable set $C\subset\{x_\alpha:\alpha<\omega_1\}$;
then $C\subset\{x_\alpha:\alpha<\gamma\}$ for some $\gamma$, and the sequence does not split at $\gamma$.

(c) The ``if'' part follows from (b); the ``only if'' part is a theorem of Arhangel'skii, see
\cite[Theorem 7.12]{Hod}.

(d) A free sequence $(x_\alpha)$ is a discrete set: $x_\alpha\notin\cl{\{x_\beta:\beta<\alpha\}}$ since the sequence
splits at $\alpha$, and $x_\alpha\notin\cl{\{x_\beta:\beta>\alpha\}}$ since it splits at $\alpha+1$.
\end{proof}

\begin{lemma}\label{l:prod}
If $K_n$, $n\in\omega$, are compact spaces of countable tightness, then $t\bigl(\prod_{n}K_n\bigr)\le\omega$.
\end{lemma}

\begin{proof}
By the theorem of Malykhin \cite{Mal}, finite products $P_n=\prod_{i<n}K_i$ have countable tightness.
Let $\pi_n\colon \prod_iK_i\to P_n$ be the projections, let $A\subset\prod_iK_i$ and $p\in\cl A$.
Then $\pi_n(p)\in\cl{\pi_n(A)}$, so there is a countable $B_n\subset A$ with $\pi_n(p)\in\cl{\pi_n(B_n)}$.
Put $B=\bigcup_nB_n$. Every basic neighborhood of $p$ has the form $\pi_n^{-1}(U)$, where $U$ is a neighborhood of
$\pi_n(p)$ in $P_n$, and it meets $B_n\subset B$. Thus, $p\in\cl B$.
\end{proof}

\subsection{Left- and right-separated sets}

A set $\{z_\alpha:\alpha<\omega_1\}$ of distinct points is \emph{left-separated} (\emph{right-separated}) if
$z_\alpha\notin\cl{\{z_\beta:\beta<\alpha\}}$ (respectively, $z_\alpha\notin\cl{\{z_\beta:\beta>\alpha\}}$) for all $\alpha$.
A space $Z$ is not hereditarily separable if and only if it contains a left-separated subset of type $\omega_1$,
and $Z$ is not hereditarily Lindel\"of if and only if it contains a right-separated subset of type $\omega_1$
(Hajnal and Juh\'asz, see \cite[Theorem 3.1]{Roi} or \cite{Juh}).

\subsection{Martin's Axiom}

We use the following well-known consequence of $\MA$: \emph{every uncountable family of conditions of a ccc
partial order contains an uncountable centered subfamily}, i.e., every ccc poset has precaliber $\aleph_1$
(Kunen and Rowbottom, see \cite[Theorem 4.2]{Wei}).
We also use the $\Delta$-system lemma.

\section{Proofs of Theorems~\ref{t:A} and~\ref{t:B}}\label{s:A}

The following lemma is the combinatorial core of the proof; it is implicit in the proof of
\cite[Theorem 19]{T93}.

\begin{lemma}\label{l:core}
$(\MA)$ Let $X$ be a space, let $\{z_\alpha:\alpha<\omega_1\}\subset X$ be a set of distinct points containing no
uncountable discrete subspaces, and let $U_\alpha$ be open sets with $z_\alpha\in U_\alpha$.
\begin{enumerate}[label=\textup{(\alph*)}]
\item If $z_\beta\notin U_\alpha$ whenever $\beta<\alpha$, then there are $m\ge1$ and finite sets
$p_\xi\subset\omega_1$, $\xi<\omega_1$, such that $|p_\xi|=m$, $p_\xi<p_\eta$ for $\xi<\eta$, and for all
$\xi<\eta$ there exist $\alpha\in p_\xi$ and $\beta\in p_\eta$ with $z_\beta\in U_\alpha$.
\item If $z_\beta\notin U_\alpha$ whenever $\beta>\alpha$, then there are $m\ge1$ and finite sets
$p_\xi\subset\omega_1$, $\xi<\omega_1$, such that $|p_\xi|=m$, $p_\xi<p_\eta$ for $\xi<\eta$, and for all
$\xi<\eta$ there exist $\alpha\in p_\eta$ and $\beta\in p_\xi$ with $z_\beta\in U_\alpha$.
\end{enumerate}
\end{lemma}

\begin{proof}
(a) Let $\mathcal P$ be the set of all finite $p\subset\omega_1$ such that $z_\beta\notin U_\alpha$ whenever
$\alpha,\beta\in p$ and $\alpha<\beta$, ordered by reverse inclusion. Note that $\mathcal P$ is closed under taking subsets,
so $p,q\in\mathcal P$ are compatible if and only if $p\cup q\in\mathcal P$.

Suppose that $\mathcal P$ is ccc. By $\MA$, the family $\{\{\alpha\}:\alpha<\omega_1\}$ contains an uncountable
centered subfamily $\{\{\alpha\}:\alpha\in A\}$; then every finite subset of $A$ belongs to $\mathcal P$.
Let $\alpha,\beta\in A$, $\beta\ne\alpha$. If $\beta<\alpha$, then $z_\beta\notin U_\alpha$ by assumption, and if
$\beta>\alpha$, then $z_\beta\notin U_\alpha$ because $\{\alpha,\beta\}\in\mathcal P$. Thus,
$U_\alpha\cap\{z_\beta:\beta\in A\}=\{z_\alpha\}$, and $\{z_\beta:\beta\in A\}$ is an uncountable discrete subspace,
a contradiction.

Hence $\mathcal P$ has an uncountable antichain. Applying the $\Delta$-system lemma and passing to an uncountable
subfamily, we obtain an antichain $\{q_\xi:\xi<\omega_1\}$ forming a $\Delta$-system with root $r$ such that the sets
$p_\xi=q_\xi\setminus r$ have the same cardinality $m$ and $p_\xi<p_\eta$ for $\xi<\eta$ (the sets $p_\xi$ are pairwise
disjoint, so every countable set of ordinals meets only countably many of them, and the required subfamily is chosen
by recursion).
Let $\xi<\eta$. Since $q_\xi\cup q_\eta\notin\mathcal P$, there are $\alpha<\beta$ in $q_\xi\cup q_\eta$ with
$z_\beta\in U_\alpha$. Since $q_\xi,q_\eta\in\mathcal P$, the points $\alpha$ and $\beta$ cannot both lie in $q_\xi$ or both
in $q_\eta$; as every element of $r$ belongs to both $q_\xi$ and $q_\eta$, neither $\alpha$ nor $\beta$ is in $r$. Thus,
one of them is in $p_\xi$ and the other is in $p_\eta$; since $\alpha<\beta$ and $p_\xi<p_\eta$, we have $\alpha\in p_\xi$
and $\beta\in p_\eta$. In particular, $m\ge1$.

(b) The proof is the same, with $\mathcal P$ consisting of all finite $p$ such that $z_\beta\notin U_\alpha$ whenever
$\alpha,\beta\in p$ and $\beta<\alpha$. For $\xi<\eta$, incompatibility gives $\beta<\alpha$ in $q_\xi\cup q_\eta$ with
$z_\beta\in U_\alpha$, and then $\beta\in p_\xi$, $\alpha\in p_\eta$.
\end{proof}

\begin{proof}[Proof of Theorem~\ref{t:A}]
We show that every $Z\subset X$ containing no uncountable discrete subspaces is hereditarily separable and
hereditarily Lindel\"of; this implies that $X$ is ESLC.

\smallskip
\emph{$Z$ is hereditarily separable.} Suppose not. Then $Z$ contains a left-separated set
$\{z_\alpha:\alpha<\omega_1\}$. Put $F_\alpha=\cl{\{z_\beta:\beta<\alpha\}}$ (the closure in $X$); then $z_\alpha\notin F_\alpha$.
Since $\Phi$ generates the topology of $X$, there are $k_\alpha\in\N$, $\bar f_\alpha\in\Phi^{k_\alpha}$ and open boxes
$I_\alpha,J_\alpha\subset\R^{k_\alpha}$ with rational vertices such that $\cl{I_\alpha}\subset J_\alpha$,
$\bar f_\alpha(z_\alpha)\in I_\alpha$ and $\bar f_\alpha(F_\alpha)\cap J_\alpha=\es$. Passing to an uncountable subset
(which is again left-separated) and reindexing, we may assume that $k_\alpha=k$, $I_\alpha=I$ and $J_\alpha=J$ for all $\alpha$.
Put $U_\alpha=\bar f_\alpha^{-1}(I)$. If $\beta<\alpha$, then $z_\beta\in F_\alpha$, hence $\bar f_\alpha(z_\beta)\notin J$ and
$z_\beta\notin U_\alpha$. By Lemma~\ref{l:core}(a), there exist $m$ and $p_\xi$, $\xi<\omega_1$, as there.

Let $p_\xi=\{\alpha^\xi_0<\alpha^\xi_1<\dots<\alpha^\xi_{m-1}\}$ and put
\[
G_\xi=\bigl(\bar f_{\alpha^\xi_0},\dots,\bar f_{\alpha^\xi_{m-1}}\bigr)\in\Phi^{km}.
\]
We write elements of $C_p(X)^{km}$ as $\bar h=(\bar h_0,\dots,\bar h_{m-1})$ with $\bar h_i\in C_p(X)^k$. For
$\xi<\omega_1$, consider the open sets
\begin{align*}
W^1_\xi&=\{\bar h\in C_p(X)^{km}: \bar h_i(z_\beta)\in I\text{ for some } i<m \text{ and } \beta\in p_\xi\},\\
W^0_\xi&=\{\bar h\in C_p(X)^{km}: \bar h_i(z_\beta)\in J\text{ for some } i<m \text{ and } \beta\in p_\xi\}.
\end{align*}
Since $W^1_\xi$ is a finite union of sets of the form $\{\bar h:\bar h_i(z_\beta)\in I\}$ and $\cl I\subset J$, we have
$\cl{W^1_\xi}\subset W^0_\xi$. We claim that
\begin{enumerate}[label=\textup{(\roman*)}]
\item $G_\eta\in W^1_\xi$ for $\eta\le\xi$;
\item $G_\eta\notin W^0_\xi$ for $\eta>\xi$.
\end{enumerate}
Indeed, $G_\xi\in W^1_\xi$, since $\bar f_\alpha(z_\alpha)\in I$ for $\alpha=\alpha^\xi_0$. If $\eta<\xi$, then there are
$\alpha\in p_\eta$ and $\beta\in p_\xi$ with $z_\beta\in U_\alpha$, i.e., $\bar f_\alpha(z_\beta)\in I$; since $\bar f_\alpha$ is
one of the blocks of $G_\eta$, we obtain $G_\eta\in W^1_\xi$. If $\eta>\xi$, then $\beta<\alpha$ for all $\alpha\in p_\eta$ and
$\beta\in p_\xi$; hence $z_\beta\in F_\alpha$ and $\bar f_\alpha(z_\beta)\notin J$, i.e., $G_\eta\notin W^0_\xi$.

By (i), $\cl{\{G_\eta:\eta\le\xi\}}\subset\cl{W^1_\xi}\subset W^0_\xi$, and by (ii),
$\cl{\{G_\eta:\eta>\xi\}}\cap W^0_\xi=\es$ (the closures are taken in $\Phi^{km}$). Thus, the sequence $(G_\xi)$ splits at
every successor ordinal, and by Lemma~\ref{l:succ}, $(G_{\xi+1})_{\xi<\omega_1}$ is a free sequence in $\Phi^{km}$.
This contradicts $\Phi^{km}\in\NFS$.

\smallskip
\emph{$Z$ is hereditarily Lindel\"of.} Suppose not. Then $Z$ contains a right-separated set $\{z_\alpha:\alpha<\omega_1\}$.
Put $F_\alpha=\cl{\{z_\beta:\beta>\alpha\}}$ and choose $k$, $\bar f_\alpha$, $I$, $J$ and $U_\alpha$ as above. Now
$z_\beta\notin U_\alpha$ for $\beta>\alpha$, and Lemma~\ref{l:core}(b) gives $m$ and $p_\xi$. Define $G_\xi$, $W^1_\xi$ and
$W^0_\xi$ by the same formulas. Then
\begin{enumerate}[label=\textup{(\roman*$'$)}]
\item $G_\eta\in W^1_\xi$ for $\eta\ge\xi$;
\item $G_\eta\notin W^0_\xi$ for $\eta<\xi$.
\end{enumerate}
Indeed, $G_\xi\in W^1_\xi$ as before. If $\eta>\xi$, there are $\alpha\in p_\eta$ and $\beta\in p_\xi$ with
$\bar f_\alpha(z_\beta)\in I$, whence $G_\eta\in W^1_\xi$. If $\eta<\xi$, then $\beta>\alpha$ for all $\alpha\in p_\eta$ and
$\beta\in p_\xi$; hence $z_\beta\in F_\alpha$ and $\bar f_\alpha(z_\beta)\notin J$, so $G_\eta\notin W^0_\xi$.

Consequently, $\cl{\{G_\eta:\eta\ge\xi\}}\subset W^0_\xi$ and $\cl{\{G_\eta:\eta<\xi\}}\cap W^0_\xi=\es$ for every $\xi$,
i.e., $(G_\xi)_{\xi<\omega_1}$ is a free sequence in $\Phi^{km}$, a contradiction.
\end{proof}

\subsection{Disjunctive families}

Todor\v{c}evi\'c's summation trick, which allows one to obtain a free sequence in $\Phi$ itself rather than in its
finite powers, works in a more general situation. What is actually used is that
a finite family of functions can be replaced by a single function of the family which realizes
the ``disjunction'' of their conditions. Let $I,J\subset\R$ be open sets with $\cl I\subset J$. We say that
$\Phi\subset C_p(X)$ is \emph{$(I,J)$-disjunctive} if for any $f_1,\dots,f_m\in\Phi$ there is $g\in\Phi$ such that
\[
\bigcup_{i\le m}f_i^{-1}(I)\subset g^{-1}(I)\qquad\text{and}\qquad g^{-1}(J)\subset\bigcup_{i\le m}f_i^{-1}(J).
\]

\begin{theorem}\label{t:Bgen}
$(\MA)$ Let $X$ be a space, let $\Phi\subset C_p(X)$, and let $\mathcal P$ be a countable family of pairs $(I,J)$
of open subsets of $\R$ with $\cl I\subset J$. Suppose that $\Phi$ is $(I,J)$-disjunctive for every $(I,J)\in\mathcal P$
and that for any closed $F\subset X$ and $x\in X\setminus F$ there are $f\in\Phi$ and $(I,J)\in\mathcal P$ with
$f(x)\in I$ and $f(F)\cap J=\es$. If $\Phi\in\NFS$, then $X$ is ESLC.
\end{theorem}

The hypotheses of Theorem~\ref{t:Bgen} hold, with $\mathcal P$ consisting of the pairs $\bigl((a,\infty),(b,\infty)\bigr)$,
$a>b$ rational, in each of the following cases:
\begin{itemize}
\item $\Phi$ is closed under finite maxima, and for any closed $F$ and $x\notin F$ there is $f\in\Phi$
with $f(x)>\sup f(F)$ (take $g=\max_i f_i$);
\item $\Phi$ consists of nonnegative functions, is closed under finite sums, and for any closed $F$ and $x\notin F$
there is $f\in\Phi$ with $f(x)>0$ and $f(F)\subset\{0\}$ (take $g=\sum_i f_i$ and $b=0$);
\item $\Phi\subset C_p(X,[0,1])$ is closed under the operation $(f,h)\mapsto f+h-fh$ and has the same separation
property (take $g=1-\prod_i(1-f_i)$ and $b=0$).
\end{itemize}
Families closed under finite minima are treated symmetrically. Theorem~\ref{t:B} is the case of sums (if $f(x)\ge1$
and $f(F)\subset\{0\}$, then $f(x)>a$ for rational $a\in(0,1)$), and
Theorem~\ref{t:Bgen} applies directly to $\Phi=C_p(X)$ with maxima, which gives another derivation of Theorem~\ref{t:T19}.

\begin{proof}[Proof of Theorem~\ref{t:Bgen}]
We repeat the proof of Theorem~\ref{t:A} with $k=1$ and the following changes. For a left-separated (right-separated)
set $\{z_\alpha\}\subset Z$ and the sets $F_\alpha$ as above, choose $f_\alpha\in\Phi$ and $(I_\alpha,J_\alpha)\in\mathcal P$
with $f_\alpha(z_\alpha)\in I_\alpha$ and $f_\alpha(F_\alpha)\cap J_\alpha=\es$. Since $\mathcal P$ is countable, passing
to an uncountable subset we may assume that $(I_\alpha,J_\alpha)=(I,J)$ for all $\alpha$. Put $U_\alpha=f_\alpha^{-1}(I)$.
The hypothesis of Lemma~\ref{l:core}(a) (respectively, (b)) is satisfied; let $m$ and $p_\xi$ be as there.
Since $\Phi$ is $(I,J)$-disjunctive, there are $g_\xi\in\Phi$ such that
\[
\bigcup_{\alpha\in p_\xi}f_\alpha^{-1}(I)\subset g_\xi^{-1}(I)\quad\text{and}\quad
g_\xi^{-1}(J)\subset\bigcup_{\alpha\in p_\xi}f_\alpha^{-1}(J).
\]
Consider the open sets
\begin{align*}
W^1_\xi&=\{h\in C_p(X): h(z_\beta)\in I\text{ for some }\beta\in p_\xi\},\\
W^0_\xi&=\{h\in C_p(X): h(z_\beta)\in J\text{ for some }\beta\in p_\xi\}.
\end{align*}
Clearly, $\cl{W^1_\xi}\subset W^0_\xi$. We have $g_\xi\in W^1_\xi$, because $f_\alpha(z_\alpha)\in I$ for $\alpha\in p_\xi$.
If there are $\alpha\in p_\eta$ and $\beta\in p_\xi$ with $z_\beta\in U_\alpha$, i.e., $f_\alpha(z_\beta)\in I$, then
$g_\eta(z_\beta)\in I$ and $g_\eta\in W^1_\xi$. If $z_\beta\in F_\alpha$, i.e., $f_\alpha(z_\beta)\notin J$, for all
$\alpha\in p_\eta$ and $\beta\in p_\xi$, then $g_\eta(z_\beta)\notin J$ for $\beta\in p_\xi$, i.e., $g_\eta\notin W^0_\xi$.
Thus, conditions (i), (ii) (respectively, (i$'$), (ii$'$)) hold for the sequence $(g_\xi)$ in $\Phi$, and we obtain a
free sequence in $\Phi$, a contradiction.
\end{proof}

\begin{remark}
Theorem~\ref{t:A} has the same structure: there the role of $g_\xi$ is played by the tuple
$(\bar f_\alpha)_{\alpha\in p_\xi}\in\Phi^{km}$, and the disjunction is realized in a finite power of $\Phi$
rather than in $\Phi$ itself. This is why Theorem~\ref{t:A} requires $\Phi^n\in\NFS$ for all $n$, whereas for
disjunctive families $\Phi\in\NFS$ suffices.
\end{remark}

\section{Proof of Theorem~\ref{t:main}}\label{s:main}

For $x\in X$, let $\hat x\colon C_p(X)\to\R$ be the evaluation function, $\hat x(f)=f(x)$. The map
$e\colon X\to C_p(C_p(X))$, $e(x)=\hat x$, is continuous, and the family $\Phi=e(X)=\{\hat x:x\in X\}$ generates
the topology of $C_p(X)$, because $C_p(X)$ carries the topology of a subspace of $\R^X$. Since $\Phi^n=e^n(X^n)$ is
a continuous image of $X^n$, Lemma~\ref{l:NFS}(a) implies that $\Phi^n\in\NFS$ for all $n$. It remains to apply
Theorem~\ref{t:A} to the space $C_p(X)$ and the family $\Phi$; thus, $C_p(X)$ is ESLC.

Finally, $C_p(X)^\omega$ is homeomorphic to $C_p(X\times\omega)$, and every finite power of $X\times\omega$ is a
countable disjoint sum of copies of a finite power of $X$. A countable disjoint sum of spaces in $\NFS$ is in $\NFS$,
since an uncountable subsequence of a free sequence lying in one summand is free. Applying what we have proved to $X\times\omega$, we see that
$C_p(X)^\omega=C_p(X\times\omega)$ is ESLC. \qed

\begin{remark}
The family $e(X)\subset C_p(C_p(X))$ is not disjunctive in general, so Theorem~\ref{t:Bgen} cannot be applied
here, and we need $\NFS$ for all finite powers of $X$. We do not know whether $X\in\NFS$ suffices
(see Question~\ref{q:1}).
\end{remark}

\section{Lindel\"of $\Sigma$-spaces}\label{s:Sigma}

Recall that a space $X$ is a \emph{$k$-space} if a set $A\subset X$ is closed whenever $A\cap C$ is closed
in $C$ for every compact $C\subset X$.

\begin{lemma}\label{l:ktight}
If $X$ is a $k$-space and all compact subspaces of $X$ have countable tightness, then $t(X)\le\omega$.
\end{lemma}

\begin{proof}
For $A\subset X$, let $[A]_\omega=\bigcup\{\cl B: B\subset A,\ |B|\le\omega\}$. It suffices to show that $[A]_\omega$ is
closed. Let $C\subset X$ be compact and $q\in C$ be in the closure of $[A]_\omega\cap C$. Since $t(C)\le\omega$, there are
$q_n\in[A]_\omega\cap C$ with $q\in\cl{\{q_n:n\in\omega\}}$. Choose countable $B_n\subset A$ with $q_n\in\cl{B_n}$;
then $q\in\cl{\bigcup_nB_n}$, i.e., $q\in[A]_\omega$. Thus, $[A]_\omega\cap C$ is closed in $C$ for every compact $C$,
and $[A]_\omega$ is closed.
\end{proof}

\begin{lemma}[{see \cite[Theorem 3.22]{Gru}}]\label{l:kperfect}
Every Lindel\"of $p$-space is a $k$-space.
\end{lemma}

\begin{proof}[Sketch of proof]
Let $\pi\colon L\to S$ be a perfect map of $L$ onto a metrizable space $S$, and let $e\colon L\to\beta L$ be the
canonical embedding. The diagonal map $\pi\triangle e\colon L\to S\times\beta L$ is perfect
\cite[Ch.~VI, Problem~55]{AP} and injective; hence it is a homeomorphism of $L$ onto a closed subspace of
$S\times\beta L$. The product of a $k$-space and a compact space is a $k$-space \cite[3.3.27]{Eng}
(see also \cite[Ch.~VI, Problem~25]{AP}), and closed subspaces of $k$-spaces are $k$-spaces \cite[Sect.~3.3]{Eng}.
\end{proof}

\begin{lemma}\label{l:graph}
Let $X$ be a Lindel\"of $\Sigma$-space. Then there exist a separable metrizable space $S$ and a subspace $M\subset S\times X$
such that the projections $\pi_S\colon M\to S$ and $\pi_X\colon M\to X$ are surjective and $\pi_S$ is perfect.
In particular, $M$ is a Lindel\"of $p$-space and every compact subset $C$ of $M$ is contained in the product
$\pi_S(C)\times\pi_X(C)$ of a metrizable compact space and a compact subspace of $X$.
\end{lemma}

\begin{proof}
Let $\pi\colon L\to S$ be a perfect map of a space $L$ onto a separable metrizable space $S$, and let $g\colon L\to X$
be a continuous surjection. Put $h=\pi\triangle g\colon L\to S\times X$ and $M=h(L)$. Then $\pi_X(M)=g(L)=X$ and
$\pi_S\circ h=\pi$. The map $\pi_S|_M$ is closed, since $\pi_S(F)=\pi(h^{-1}(F))$ for $F\subset M$, and its fibers
$\pi_S^{-1}(s)\cap M=h(\pi^{-1}(s))$ are compact. Thus, $\pi_S|_M$ is perfect. (In fact, $h$ itself is perfect
\cite[Ch.~VI, Problem~55]{AP}, so $M$ is closed in $S\times X$.)
\end{proof}

\begin{proof}[Proof of Theorem~\ref{t:Sigma}]
(1)$\Rightarrow$(2). A compact subspace $K\subset X$ is closed, so $K\in\NFS$ by Lemma~\ref{l:NFS}(a), and $t(K)\le\omega$
by Lemma~\ref{l:NFS}(c).

(2)$\Rightarrow$(3). Let $M$ be as in Lemma~\ref{l:graph}; it is a Lindel\"of $p$-space mapped continuously onto $X$.
By Lemma~\ref{l:kperfect}, $M$ is a $k$-space. If $C\subset M$ is compact, then $C$ is a subspace of
$\pi_S(C)\times\pi_X(C)$, which has countable tightness by (2) and Lemma~\ref{l:prod}; hence $t(C)\le\omega$.
By Lemma~\ref{l:ktight}, $t(M)\le\omega$.

(3)$\Rightarrow$(1). A Lindel\"of $p$-space of countable tightness is in $\NFS$ by Lemma~\ref{l:NFS}(b), and so is its
continuous image $X$ by Lemma~\ref{l:NFS}(a).

(2)$\Rightarrow$(4). The space $X^\omega$ is Lindel\"of $\Sigma$ \cite{Arh}. If $C\subset X^\omega$ is compact, then $C$ is a subspace
of $\prod_n\pi_n(C)$, where $\pi_n\colon X^\omega\to X$ are the projections; by (2) and Lemma~\ref{l:prod}, $t(C)\le\omega$.
Thus, $X^\omega$ satisfies (2), and by the implication (2)$\Rightarrow$(1) for $X^\omega$, we have $X^\omega\in\NFS$.

(4)$\Rightarrow$(1). $X$ is a continuous image of $X^\omega$.
\end{proof}

\section{Proofs of the corollaries}\label{s:cor}

\begin{proof}[Proof of Corollary~\ref{c:Sigma}]
Let $Y\subset C_p(X)$ contain no uncountable discrete subspaces. By Theorem~\ref{t:Sigma}, $X^\omega\in\NFS$, so all finite powers of $X$ are in $\NFS$. By Theorem~\ref{t:main}, $Y$ is
hereditarily separable. By Proposition~\ref{p:known}(1), $Y$ has a countable network.
\end{proof}

\begin{proof}[Proof of Corollary~\ref{c:L}]
Let $Y\subset C_p(X)$ contain no uncountable discrete subspaces.

(1) By Lemma~\ref{l:NFS}(b), $X\in\NFS$, and it remains to apply Corollary~\ref{c:Sigma}.

(2) By Proposition~\ref{p:known}(3), $Y$ is hereditarily Lindel\"of; in particular, $Y$ is Lindel\"of. Consider the
map $\psi\colon X\to C_p(Y)$, $\psi(x)(y)=y(x)$, and put $Z=\psi(X)$. The space $Z$ is a Lindel\"of $\Sigma$-space,
being a continuous image of $X$. Every compact subspace of $Z$ is a compact subspace of $C_p(Y)$ for a Lindel\"of
space $Y$, and under $\PFA$ every such compact space has countable tightness \cite[Theorem IV.11.14]{Arh}. By Theorem~\ref{t:Sigma},
$Z\in\NFS$.

For $y\in Y$, the function $j(y)\colon Z\to\R$, $j(y)(\psi(x))=y(x)$, is well defined and continuous, being the
restriction to $Z$ of the evaluation function at $y$ on $C_p(Y)$. The map $j\colon Y\to C_p(Z)$ is a homeomorphic
embedding: it is injective, and both topologies on $Y$ are generated by the functions $y\mapsto y(x)=j(y)(\psi(x))$,
$x\in X$. Since $\PFA$ implies $\MA$, Corollary~\ref{c:Sigma} applied to $Z$ and $j(Y)\subset C_p(Z)$ shows that $Y$
has a countable network.
\end{proof}

We show that Okunev's theorem \cite[Theorem 1]{Ok} (proved there under $\mathrm{MA}+\neg\mathrm{CH}$)
splits into a ``descriptive'' part, which is Theorem 1.1 of \cite{OR}, and a ``combinatorial'' part, which is covered
by Theorem~\ref{t:main}.

\begin{corollary}[Okunev]\label{c:Ok}
$(\MA)$ Let $K$ be a separable compact space, and let $Y\subset C_p(K)$ be such that $Y^n$ is Lindel\"of for all $n$.
Then $Y$ has a countable network.
\end{corollary}

\begin{proof}
Let $\psi\colon K\to C_p(Y)$, $\psi(x)(y)=y(x)$, and $K_1=\psi(K)$. Then $K_1$ is a separable compact space, and
$Y$ is homeomorphic to a subspace $Y_1$ of $C_p(K_1)$ (see \cite[p.~92]{Ok} or the proof of Corollary~\ref{c:L}).
Since all finite powers of $Y$ are Lindel\"of, $t(C_p(Y))\le\omega$ by the Arhangel'skii--Pytkeev theorem
\cite[Theorem II.1.1]{Arh}; hence $t(K_1)\le\omega$. By \cite[Theorem 1.1]{OR}, the Lindel\"of subspace $Y_1$ of
$C_p(K_1)$ is hereditarily Lindel\"of; in particular, it contains no uncountable discrete subspaces.
By Corollary~\ref{c:L}(1), $Y_1$ has a countable network.
\end{proof}

\section{Corson and Gul'ko compact spaces}\label{s:Corson}

Recall that a compact space $X$ is \emph{Corson compact} if it embeds in the $\Sigma$-product
$\Sigma(\R^A)=\{x\in\R^A:|\{a:x(a)\ne0\}|\le\omega\}$ for some $A$, and \emph{Gul'ko compact} if $C_p(X)$ is a
Lindel\"of $\Sigma$-space; every Gul'ko compact space is Corson compact. Let $C_{p,0}(X)=X$ and
$C_{p,n+1}(X)=C_p(C_{p,n}(X))$.

\begin{theorem}\label{t:Corson}
$(\MA)$ Let $X$ be a Corson compact space. Then for every $n\in\omega$ all finite and countable powers of
$C_{p,n}(X)$ are ESLC.
\end{theorem}

\begin{proof}
By Sokolov's theorem \cite[Corollary 2]{Sok}, all spaces $C_{p,k}(X)$, $k\ge1$, are Lindel\"of. Let $Z=C_{p,n}(X)$.
Since $C_p(Z)$ and $C_p(C_p(Z))$ are Lindel\"of, Asanov's theorem \cite[Theorem I.4.1]{Arh} gives $t(Z^m)\le\omega$
and $t(C_p(Z))\le\omega$; the latter means that all $Z^m$ are Lindel\"of by the Arhangel'skii--Pytkeev theorem
\cite[Theorem II.1.1]{Arh}. Hence $Z^m\in\NFS$ for all $m$ (Lemma~\ref{l:NFS}(b)), and it remains to apply
Todor\v{c}evi\'c's theorem in the form stated in the Introduction.
\end{proof}

For $n=1$, Theorem~\ref{t:Corson} is contained in Corollary~\ref{c:L}(1); for $n\ge2$ it is not covered by Corollary~\ref{c:L}, since $C_p(X)$ is not
Lindel\"of $\Sigma$ unless $X$ is Gul'ko compact. The proof works for every space $X$ such that all $C_{p,k}(X)$,
$k\ge1$, are Lindel\"of, e.g., for every Sokolov space $X$ with all finite powers Lindel\"of and countably tight
\cite{Tk05}, \cite[Problem 162]{Tk15}.

\begin{proposition}\label{p:CorsonZFC}
\begin{enumerate}[label=\textup{(\arabic*)}]
\item If $X$ is Corson compact, then no $C_{p,n}(X)$, $n\in\omega$, contains an $S$-space.
\item If $X$ is Gul'ko compact, then $X$ is ESLC.
\end{enumerate}
\end{proposition}

\begin{proof}
(1) The space $X$ is $\aleph_0$-monolithic, since the closure of a countable subset of $\Sigma(\R^A)$ lies in a
countable subproduct, and $\aleph_0$-stable, being compact. By the duality theorems recalled in
Section~\ref{s:prel}, all $C_{p,n}(X)$ are $\aleph_0$-monolithic and $\aleph_0$-stable. In an $\aleph_0$-monolithic
space every hereditarily separable subspace has a countable network and is therefore hereditarily Lindel\"of.

(2) If $Y\subset X$ is hereditarily separable or hereditarily Lindel\"of, then the Gul'ko compact space $\cl Y$ has the
Souslin property and is therefore metrizable \cite[Problem 294]{Tk15}.
\end{proof}

Let $L(\kappa)=\kappa\cup\{\infty\}$ denote the one-point Lindel\"ofication of the discrete space $\kappa$: the points of
$\kappa$ are isolated, and the neighborhoods of $\infty$ are the sets containing $\infty$ with countable complement.

\begin{lemma}\label{l:Lkappa}
Let $\varphi\colon L(\kappa)^n\to Y$ be a continuous map to a Hausdorff space $Y$, and let $Z\subset Y$ be a Lindel\"of
subspace of countable pseudocharacter. Then $Z\cap\varphi(L(\kappa)^n)$ is countable. In particular, every
hereditarily Lindel\"of subspace of a Hausdorff continuous image of $L(\kappa)^n$ is countable.
\end{lemma}

\begin{proof}
We use induction on $n$; for $n=0$ the set $\varphi(L(\kappa)^0)$ is a point. Let $n\ge1$ and $p=\varphi(\infty,\dots,\infty)$.
Every $z\in Z\setminus\{p\}$ has a neighborhood $U_z$ disjoint from a neighborhood $V_z$ of $p$. The set
$\varphi^{-1}(V_z)$ contains a basic neighborhood $(\{\infty\}\cup(\kappa\setminus C_z))^n$ of $(\infty,\dots,\infty)$,
where $C_z\subset\kappa$ is countable; hence $\varphi^{-1}(U_z)\subset\bigcup_{i<n}\{x:x_i\in C_z\}$. Since $Z$ has
countable pseudocharacter, $Z\setminus\{p\}$ is an $F_\sigma$-subset of the Lindel\"of space $Z$, so it is Lindel\"of and
is covered by countably many sets $U_z$. Let $C$ be the union of the corresponding sets $C_z$. Then
\[
\bigl(Z\cap\varphi(L(\kappa)^n)\bigr)\setminus\{p\}\subset\bigcup_{i<n,\ c\in C}\varphi(S_{i,c}),\qquad
S_{i,c}=\{x\in L(\kappa)^n:x_i=c\}.
\]
Each $S_{i,c}$ is homeomorphic to $L(\kappa)^{n-1}$, so $Z\cap\varphi(S_{i,c})$ is countable by the inductive
assumption applied to $\varphi|_{S_{i,c}}$. Thus, $Z\cap\varphi(L(\kappa)^n)$ is countable.
\end{proof}

\begin{proposition}\label{p:CpD}
Let $X$ be a zero-dimensional Corson compact space and $D=\{0,1\}$. Then every Lindel\"of subspace of $C_p(X,D)$ of
countable pseudocharacter is countable. In particular, every hereditarily Lindel\"of subspace of $C_p(X,D)$ is countable,
and $C_p(X,D)$ is ESLC in $\mathrm{ZFC}$.
\end{proposition}

\begin{proof}
The space $X$ has a point-countable $T_0$-separating family $\mathcal C$ of clopen sets (see \cite[Problems 118,
150]{Tk15}). Every neighborhood of $0$ in $C_p(X,D)$ contains all but countably many characteristic functions $\chi_C$,
$C\in\mathcal C$; hence $\{\chi_C:C\in\mathcal C\}\cup\{0\}$ is a continuous image of $L(\kappa)$, $\kappa=|\mathcal C|$.
Every clopen subset of $X$ is a Boolean combination of finitely many elements of $\mathcal C$, and Boolean operations are
continuous on $C_p(X,D)$. Thus, $C_p(X,D)$ is a countable union of continuous images of finite powers of $L(\kappa)$
(cf.\ \cite[solution of Problem 150]{Tk15}). If $Z\subset C_p(X,D)$ is Lindel\"of and has countable pseudocharacter, then
by Lemma~\ref{l:Lkappa} its intersection with each of these images is countable, so $Z$ is countable; a hereditarily
Lindel\"of subspace of the Hausdorff space $C_p(X,D)$ has countable pseudocharacter. Finally, a hereditarily separable subspace of
$C_p(X,D)\subset C_p(X)$ is hereditarily Lindel\"of by Proposition~\ref{p:CorsonZFC}(1).
\end{proof}

The argument does not apply to $C_p(X)$. The space $C_p(X)$ is a continuous image of a closed subspace of
$L(\kappa)^\omega$ for every Corson compact space $X$ \cite[Problem 150]{Tk15}, but the proof of Lemma~\ref{l:Lkappa}
breaks down for $L(\kappa)^\omega$: a basic neighborhood of $(\infty,\infty,\dots)$ restricts only finitely many
coordinates, and the sets $\{x:x_i=c\}$ are again homeomorphic to $L(\kappa)^\omega$.
Of course, for $L(\kappa)^\omega$ the conclusion of Lemma~\ref{l:Lkappa} must be weakened, since $[0,1]$ is a continuous
image of $L(\kappa)^\omega$ ($L(\kappa)^\omega$ is normal and contains a closed copy of $\omega^\omega$); the expected
analogue is stated in Question~\ref{q:Lomega}.

\begin{proposition}\label{p:CorsonL}
Suppose that there is a non-separable perfectly normal compact space \textup(e.g., the Souslin continuum, if there is
a Souslin line, or a space constructed under $\diamondsuit$ \cite[Problem 131]{Tk15}\textup). Then there is a
perfectly normal non-metrizable Corson compact space $X$; it is an $L$-space, and $C_{p,2n}(X)$ is not ESLC for every
$n\in\omega$.
\end{proposition}

\begin{proof}
A non-separable perfectly normal compact space $K$ is first countable, so by Shapirovski\u{\i}'s theorem it admits an
irreducible map $f$ onto a Corson compact space $X$ \cite[Problems 129, 130]{Tk15}. Then $X$ is hereditarily
Lindel\"of and $d(X)=d(K)>\omega$, i.e., $X$ is an $L$-space. Finally, $X$ embeds in $C_{p,2n}(X)$ by evaluation.
\end{proof}

For such $X$, the space $C_p(X)$ contains no $S$-spaces (Proposition~\ref{p:CorsonZFC}(1)), but has uncountable spread
\cite[Problem 128]{Tk15}. Under $\MA$ such spaces $X$ do not exist, since under $\MA$ a Corson compact space with the
Souslin property has precaliber $\aleph_1$ and hence is metrizable \cite[Problem 132]{Tk15}. Thus, in Theorem~\ref{t:Corson} the assumption $\MA$ cannot be
omitted for even $n$; the case of odd $n$ is open (Question~\ref{q:Corson}).

\section{What happens without $\MA$?}\label{s:noMA}

All main results of this paper except Theorem~\ref{t:Sigma} use $\MA$, and this assumption cannot be omitted. The
counterexamples are based on the following result of Todor\v{c}evi\'c \cite{T89}: if $\mathfrak b=\omega_1$
(in particular, under $\mathrm{CH}$), then there exists a non-metrizable compact space $K$ all of whose finite powers
are hereditarily separable. Under $\mathrm{CH}$, one can take for $K$ the one-point compactification of the
Kunen line \cite{JKR}. In the terminology of \cite{HK}, $K^2$ is then a \emph{strong $S$-space}, i.e., all its finite
powers are hereditarily separable, while $K^2$ is not hereditarily Lindel\"of (see the proof below); recall that
strong $S$- and $L$-spaces do not exist under $\MA$ \cite{Kun}.

\begin{proposition}\label{p:noMA}
Let $K$ be a non-metrizable compact space all of whose finite powers are hereditarily separable.
\begin{enumerate}[label=\textup{(\arabic*)}]
\item Let $X=K^2$ and $\Phi=C_p(K^2)$. Then $\Phi$ generates the topology of $X$, $\Phi^n\in\NFS$ for all $n$, and
$X$ is hereditarily separable but not hereditarily Lindel\"of. Thus, Theorems~\ref{t:T19}, \ref{t:A}, \ref{t:B}
and~\ref{t:Bgen} fail for $X$ and $\Phi$ \textup(for Theorem~\ref{t:B}, with $\Phi=C_p^+(K^2)$\textup).
\item Let $X=C_p(K^2)$. Then $X^n\in\NFS$ for all $n$, and $C_p(X)$ contains the subspace $K^2$, which is hereditarily
separable but not hereditarily Lindel\"of. Thus, Theorem~\ref{t:main} fails for $X$.
\item Let $X=K$ \textup(or $X=K^2$\textup). Then $X$ is a compact space of countable tightness, so $X^n\in\NFS$ for
all $n$, and $C_p(X)$ is an $L$-space: it is hereditarily Lindel\"of but not hereditarily separable; in particular, it has no countable
network. Thus, Theorem~\ref{t:main} fails for $X$ also in the ``$L$-direction'', and Corollaries~\ref{c:Sigma}
and~\ref{c:L}\textup{(1)} fail for $X$.
\end{enumerate}
\end{proposition}

\begin{proof}
All finite powers of $K^2$ are hereditarily separable, so by the Zenor--Velichko theorem \cite{Zen, Vel} all finite
powers of $C_p(K^2)$ are hereditarily Lindel\"of, hence belong to $\NFS$ by Lemma~\ref{l:NFS}(d). On the other hand, $K^2$ is not hereditarily
Lindel\"of: otherwise the diagonal of $K$ would be a $G_\delta$-set in $K^2$, and $K$ would be metrizable by
\v{S}ne\u{\i}der's theorem \cite{Eng}.

(1) Clearly, $C_p(K^2)$ generates the topology of $K^2$, and $C_p^+(K^2)$ is closed in $C_p(K^2)$, hence belongs
to $\NFS$ together with $C_p(K^2)$.

(2) The space $K^2$ embeds in $C_p(C_p(K^2))$ via $x\mapsto\hat x$.

(3) All finite powers of $X$ are hereditarily separable, so they belong to $\NFS$ by Lemma~\ref{l:NFS}(d); in
particular, $t(X)\le\omega$ by Lemma~\ref{l:NFS}(c). The space $C_p(X)$ is hereditarily Lindel\"of by the Zenor--Velichko theorem.
It is not hereditarily separable: by Velichko's theorem \cite{Vel} (see also \cite[p.~608]{Arh90}), if $C_p(X)$ is
hereditarily separable, then all finite powers of $X$ are hereditarily Lindel\"of, while $X^2$ is not. Finally, $\mathrm{nw}(C_p(X))=\mathrm{nw}(X)=w(X)>\omega$ \cite[Theorem I.1.3]{Arh}.
\end{proof}

Okunev's theorem (Corollary~\ref{c:Ok}) also fails without additional axioms. As noted in \cite[p.~96]{Ok}, under
$\mathrm{CH}$ Pol \cite{Pol} constructed a separable scattered compact space $K$ with an uncountable discrete subspace
such that all finite powers of $Y=C_p(K,\{0,1\})$ are Lindel\"of; the space $Y\subset C_p(K)$ separates the points
of $K$, and it has no countable network, since otherwise $K\subset C_p(Y)$ would be metrizable.

\section{Questions}\label{s:q}

By Theorem~\ref{t:Bgen}, for disjunctive families (e.g., closed under finite maxima or sums) one can avoid finite
powers of $\Phi$ in Theorem~\ref{t:A}.
The general case is related to the problem of pulling down free sequences from finite powers, discussed in \cite[\S\S 8, 10, 11]{T93}.

\begin{question}\label{q:1}
$(\MA)$ Let $\Phi\subset C_p(X)$ generate the topology of $X$ and $\Phi\in\NFS$. Is it true that every subspace of $X$
with no uncountable discrete subspaces is hereditarily Lindel\"of and hereditarily separable? In particular, does
Theorem~\ref{t:main} hold under the assumption $X\in\NFS$?
\end{question}

\begin{question}[Arhangel'skii]\label{q:tight}
Let $X$ be a Lindel\"of space and let $K$ be a compact subspace of $C_p(X)$. Is it true that $t(K)\le\omega$?
\end{question}

This is Problem IV.11.12$'$ of \cite{Arh} (equivalent to Problem IV.11.11 there) and Problem 28 of \cite{Arh90}
(see also \cite[Problem 4.5.1]{Tk15}). The answer is positive under $\PFA$ \cite[Theorem IV.11.14]{Arh}; in $\mathrm{ZFC}$ it is known that $\omega_1+1$ does not embed in $C_p(X)$
for a Lindel\"of space $X$ \cite{AU}. We do not know the answer under $\MA$ (cf.\ \cite[Problem 4.5.2]{Tk15}).

In Question~\ref{q:tight} one may assume that $X$ is an Eberlein--Grothendieck space (i.e., a subspace of $C_p(K)$ for
a compact space $K$) whose topology is generated by $K$.
Indeed, let $\psi\colon X\to C_p(K)$, $\psi(x)(f)=f(x)$, and $Z=\psi(X)$. Every $f\in K$ factors through $\psi$, so
$K$ embeds in $C_p(Z)$ and generates the topology of $Z$; moreover, $Z$ is Lindel\"of if $X$ is. Such a space $Z$ has
countable tightness, since $t(C_p(K))\le\omega$ for compact $K$ by the Arhangel'skii--Pytkeev theorem
\cite[Theorem II.1.1]{Arh}, and it is $\aleph_0$-monolithic (Proposition~\ref{p:known}(1)). For subspaces $Z$ of
$C_p(K)$ with $K$ compact, the following conditions are equivalent: $Z$ is Lindel\"of; $Z\in\NFS$; $Z$ has countable
extent. Indeed, a Lindel\"of space of countable tightness is in $\NFS$ (Lemma~\ref{l:NFS}(b)); every space in $\NFS$ has
countable extent, since any enumeration of an uncountable closed discrete set of type $\omega_1$ is a free sequence;
and a subspace of $C_p(K)$ of countable extent is Lindel\"of by Baturov's theorem (Proposition~\ref{p:known}(2)).
Since $Z\in\NFS$ whenever $X\in\NFS$ (Lemma~\ref{l:NFS}(a)), and $\mathrm{ext}(Z)\le\mathrm{ext}(X)$ (choosing one
preimage of each point of a closed discrete subset of $Z$ gives a closed discrete subset of $X$), in
Question~\ref{q:tight} the assumption that $X$ is Lindel\"of may be replaced by $X\in\NFS$ or by
$\mathrm{ext}(X)\le\omega$.

\begin{question}\label{q:L}
$(\MA)$ Is it true that $C_p(X)$ contains no $L$-spaces
\begin{enumerate}[label=\textup{(\arabic*)}]
\item for every Lindel\"of $\Sigma$-space $X$?
\item for every compact space $X$?
\end{enumerate}
\end{question}

Question~\ref{q:L}(2) is closely related to Question 0.5 of \cite{OR}. If the answer to Question~\ref{q:tight} is
positive under $\MA$, then so is the answer to Question~\ref{q:L} (both parts). Indeed, let $X$ be a Lindel\"of
$\Sigma$-space and let $Y\subset C_p(X)$ be hereditarily Lindel\"of. Repeating the proof of Corollary~\ref{c:L}(2),
with Question~\ref{q:tight} used instead of $\PFA$, we see that $Z=\psi(X)\subset C_p(Y)$ is a Lindel\"of
$\Sigma$-space all of whose compact subspaces have countable tightness; hence $Z\in\NFS$ by Theorem~\ref{t:Sigma}, and
$Y$, which embeds in $C_p(Z)$, has a countable network by Corollary~\ref{c:Sigma}.

The following questions concern Section~\ref{s:Corson}.

\begin{question}\label{q:Corson}
\begin{enumerate}[label=\textup{(\arabic*)}]
\item Let $X$ be a perfectly normal non-metrizable Corson compact space. Is $C_{p,2n+1}(X)$ ESLC? In particular, does
$C_p(X)$ contain an $L$-space?
\item Is $C_{p,2n+1}(X)$ ESLC for every Corson compact space $X$?
\item Is $C_{p,n}(X)$ ESLC for every Gul'ko compact space $X$ and every $n\in\N$?
\end{enumerate}
\end{question}

By Proposition~\ref{p:CorsonZFC}(1), these spaces contain no $S$-spaces, so the question is whether they contain
$L$-spaces; for zero-dimensional $X$, the analogue of (1) for $C_p(X,D)$ has a positive answer by
Proposition~\ref{p:CpD}. Under $\MA$ the answers are positive by Theorem~\ref{t:Corson}, so Question~\ref{q:Corson} concerns
$\mathrm{ZFC}$ or models of $\neg\MA$, e.g., $\mathrm{CH}$; by Proposition~\ref{p:CorsonL}, the analogue of (2) for
even indices has a consistent negative answer.

\begin{question}\label{q:Lomega}
Let $Y$ be a continuous image of $L(\kappa)^\omega$, and let $Z\subset Y$ be hereditarily Lindel\"of.
\begin{enumerate}[label=\textup{(\arabic*)}]
\item Is $Z$ separable?
\item Does $Z$ have a countable network?
\item The same questions if $Y$ is an Eberlein--Grothendieck space.
\end{enumerate}
\end{question}

Lemma~\ref{l:Lkappa} is a $\mathrm{ZFC}$ version of Question~\ref{q:Lomega} for finite powers of $L(\kappa)$, with a
stronger conclusion. Questions (1) and (2) do not change if continuous images of closed subspaces of $L(\kappa)^\omega$
are allowed: if $Y=\varphi(F)$ with $F$ closed in $L(\kappa)^\omega$, then, since $L(\kappa)^\omega$ is normal,
$\varphi$ extends to a continuous map of $L(\kappa)^\omega$ into a Tychonoff cube containing $Y$.

The same is true for (3). Let $Y\subset C_p(K)$, $K$ compact, be a continuous image of a closed subspace of
$L(\kappa)^\omega$, and let $\psi\colon K\to C_p(Y)$, $\psi(x)(y)=y(x)$. Then $K'=\psi(K)$ is Corson compact
\cite[Problem 150]{Tk15} and $Y$ embeds in $C_p(K')$. Let $X'$ be a zero-dimensional Corson compact space which maps
onto $K'$ \cite[Problem 178]{Tk15}; then $C_p(K')$ embeds in $C_p(X')$. Finally, $C_p(X')$ is a continuous image of
$L(\kappa)^\omega$. Indeed, continuous images of $L(\kappa)^\omega$ are preserved by countable unions and countable
products (as $L(\kappa)^\omega$ maps onto $L(\kappa)^\omega\times\omega$); so this follows from the proof of
Proposition~\ref{p:CpD}, the equality $C_p(X',D)^\omega=C_p(X',D^\omega)$, and the fact that $C_p(X',[0,1])$ is a
continuous image of $C_p(X',D^\omega)$ \cite[Problem 004]{Tk15}. Thus, $Y$ is a subspace of the Eberlein--Grothendieck
space $C_p(X')$, which is a continuous image of $L(\kappa)^\omega$.

Moreover, all finite powers of $C_p(X')$ are Lindel\"of \cite[Problem 141]{Tk15} and have countable tightness, and
separable subspaces of $C_p(X')$ have a countable network, $C_p(X')$ being $\aleph_0$-monolithic. Hence, in (3) the two
questions are equivalent, and under $\MA$ the answer is positive by Todor\v{c}evi\'c's theorem. Since $C_p(X)$ is
Eberlein--Grothendieck and is a continuous image of a closed subspace of $L(\kappa)^\omega$ for every Corson compact space
$X$ \cite[Problem 150]{Tk15}, a positive answer to (3) is equivalent to the absence of $L$-spaces in $C_p(X)$ for all
Corson compact $X$, i.e., to the case $n=0$ of Question~\ref{q:Corson}(2).

Hart and Kunen \cite{HK} showed that the product of two ESLC spaces can contain both an $S$-space and an $L$-space
(under their Chain of Models Axiom, which contradicts $\MA$). Note that $C_p(X)\times C_p(Y)$ is homeomorphic to
$C_p(X\oplus Y)$, and the finite powers of $X\oplus Y$ are finite disjoint sums of spaces $X^k\times Y^l$. Hence, under
$\MA$, Theorem~\ref{t:main} shows that $C_p(X)\times C_p(Y)$ is ESLC whenever $X^k\times Y^l\in\NFS$ for all $k,l$.
Under the weaker assumption that $X^n,Y^n\in\NFS$ for all $n$, both factors are ESLC by Theorem~\ref{t:main}.

\begin{question}\label{q:prod}
$(\MA)$ Let $X$ and $Y$ be spaces such that $X^n\in\NFS$ and $Y^n\in\NFS$ for all $n$. Is $C_p(X)\times C_p(Y)$ ESLC?
\end{question}

Natural test cases are provided by Michael-type lines. Following \cite{OT}, for $A\subset X$ let $X_{(A)}$ denote the
space obtained from $X$ by retaining the topology at the points of $A$ and declaring the points of $X\setminus A$
isolated, and call $A$ \emph{holding} if $(X_{(A)})^\omega$ is Lindel\"of. Okunev and Tamano \cite[Theorem 1.5]{OT}
proved in $\mathrm{ZFC}$ that $\R$ contains $\mathfrak c$ pairwise disjoint holding sets. The diagonal of
$\R_{(A)}\times\R_{(B)}$ is closed and homeomorphic to $\R_{(A\cap B)}$ \cite[\S1]{OT}. Hence, if $A\cap B$ is
holding, then all finite products of copies of $\R_{(A)}$ and $\R_{(B)}$ are Lindel\"of and first countable, and under
$\MA$ the space $C_p(\R_{(A)})\times C_p(\R_{(B)})$ is ESLC by Theorem~\ref{t:main}. If $A\cap B=\es$, then the
diagonal is an uncountable closed discrete subset of $\R_{(A)}\times\R_{(B)}$, so this product is not in $\NFS$, and
Theorem~\ref{t:main} does not apply to $\R_{(A)}\oplus\R_{(B)}$.

\begin{question}\label{q:Michael}
Let $A,B\subset\R$ be disjoint holding sets.
\begin{enumerate}[label=\textup{(\arabic*)}]
\item Is $C_p(\R_{(A)})$ ESLC?
\item Is $C_p(\R_{(A)})\times C_p(\R_{(B)})$ ESLC?
\end{enumerate}
\end{question}

Under $\MA$, the answer to Question~\ref{q:Michael}(1) is positive by Theorem~\ref{t:main}: the finite powers of
$\R_{(A)}$ are Lindel\"of and first countable, hence belong to $\NFS$. We do not know whether $C_p(\R_{(A)})$ is
ESLC in $\mathrm{ZFC}$, and we do not know the answer to Question~\ref{q:Michael}(2) even under $\MA$.

A closely related example is \cite[Theorem 2.2]{OT}: for disjoint holding sets $A,B\subset I=[0,1]$ there are separable
scattered $\sigma$-compact spaces $X_A$, $X_B$ such that $C_p(X_A)^\omega$ and $C_p(X_B)^\omega$ are Lindel\"of, while
$C_p(X_A)\times C_p(X_B)$ contains a closed copy of $I_{(A)}\times I_{(B)}$ and is not Lindel\"of. Here $X_A$ embeds in
$C_p(M_A)$ for a space $M_A$ with Lindel\"of $M_A^\omega$, so all finite powers of $X_A$ are $\sigma$-compact and have
countable tightness by the Arhangel'skii--Pytkeev theorem; hence $X_A^n\in\NFS$, and similarly for $X_B$. Thus, $X_A$ and
$X_B$ satisfy the hypotheses of Question~\ref{q:prod}, and it would be interesting to know whether
$C_p(X_A)\times C_p(X_B)$ is ESLC under $\MA$.

Under $\MA$, every separable compact space $K$ with Lindel\"of $C_p(K)$ is metrizable, so $C_p(K)$ has a countable network.
Under $\mathrm{CH}$, the one-point compactification $K$ of the Kunen line \cite{JKR} is a non-metrizable separable compact space
all of whose finite powers are hereditarily separable; hence $C_p(K)$ is hereditarily Lindel\"of by the
Zenor--Velichko theorem \cite{Zen, Vel}. This leaves the following questions.

\begin{question}\label{q:sep}
Let $K$ be a separable compact space such that $C_p(K)$ is Lindel\"of.
\begin{enumerate}[label=\textup{(\arabic*)}]
\item Is it true that $K$ contains no uncountable discrete subspaces?
\item Is $K$ hereditarily separable?
\item Is $C_p(K)$ hereditarily Lindel\"of?
\end{enumerate}
\end{question}

These questions are related as follows. Clearly, (2) implies (1). Further, (3) implies (2), because
$\hd(X)\le\hL(C_p(X))$ for every space $X$: if $\{x_\alpha:\alpha<\omega_1\}\subset X$ is left-separated, choose
$f_\alpha\in C_p(X)$ with $f_\alpha(x_\alpha)=1$ and $f_\alpha\bigl(\cl{\{x_\beta:\beta<\alpha\}}\bigr)\subset\{0\}$;
then $f_\beta(x_\alpha)=0$ for $\beta>\alpha$, so $\{f_\alpha:\alpha<\omega_1\}$ is right-separated in $C_p(X)$.
(In the same way, $s(X)\le s(C_p(X))$.) By Baturov's theorem (Proposition~\ref{p:known}(3)), (3) is equivalent to
the absence of uncountable discrete subspaces in $C_p(K)$; thus (1) and (3) are the same question for $K$ and for
$C_p(K)$. If all finite powers of $K$ are hereditarily separable, then (3) holds by the Zenor--Velichko theorem;
whether the converse holds is a particular case of Velichko's problem \cite[Problem 18]{Arh90}.
Under $\MA$ the answer to all three questions is positive, since $K$ is metrizable. Compactness of $K$ is essential
here: Okunev and Tamano \cite[Theorem 2.1]{OT} constructed in $\mathrm{ZFC}$ a separable scattered $\sigma$-compact space
$X$ with $X^{(3)}=\es$ such that $C_p(X)^\omega$ is Lindel\"of, but $X$ contains an uncountable discrete subspace.

Question (2) is the particular case $Y=C_p(K)$ of the following question of Arhangel'skii \cite{Arh98}
(see \cite[Problem 4.5.5]{Tk15} and \cite[p.~91]{Ok}): must every separable compact subspace of $C_p(Y)$, where $Y$
is Lindel\"of, be hereditarily separable? This question is independent of $\mathrm{ZFC}$: under $\PFA$ every such
compact space is metrizable \cite[Theorem 1.8]{OR}, while under $\mathrm{CH}$ Pol \cite{Pol} constructed a separable
scattered compact space $K$ with an uncountable discrete subspace such that all finite powers of $Y=C_p(K,\{0,1\})$
are Lindel\"of; clearly, $K$ embeds in $C_p(Y)$ (see \cite[p.~96]{Ok}). However, in Pol's example the space $C_p(K)$
is not Lindel\"of \cite{OR}, so it does not answer even Question~\ref{q:sep}(1). Related questions on hereditary
Lindel\"ofness and spread of finite powers of $C_p(X)$ are Problems 18 and 19 of \cite{Arh90} (Problems 1043 and
1044 in the global numbering of the volume).

\subsection*{Declaration of generative AI use}
During the preparation of this work the author used Claude (Anthropic) in order to draft and edit the English text,
organize the exposition, search the literature, and check references and arguments. After using this tool, the author
reviewed and edited the content as needed and takes full responsibility for the content of the paper.

\end{document}